\documentclass[11pt,a4paper]{article}
\usepackage[cp1251]{inputenc}
\usepackage{amsfonts}
\usepackage{amssymb}
\usepackage{amsthm}
\usepackage{amsmath}
\usepackage{mathtext}
\usepackage{longtable}
\usepackage[russian]{babel}

\theoremstyle{remark}

\newcommand{\R}{\mathbb R}
\newcommand{\N}{\mathbb N}
\renewcommand{\P}{{\sf P}}
\newcommand{\E}{{\sf E}}
\newcommand{\D}{{\sf D}}

\renewcommand{\le}{\leqslant}
\renewcommand{\ge}{\geqslant}
\newcommand{\I}{\mathbb{I}}

\renewcommand{\phi}{\varphi}

\newcommand{\eqd}{\stackrel{d}{=}}
\newcommand{\abs}{\text{\footnotesize$|$}}
\newcommand{\vp}{\emph{\textbf{p}}}
\newcommand{\wX}{\widetilde X}

\title{Bounds of the accuracy of the normal approximation to the
distributions of random sums under relaxed moment
conditions\thanks{Research supported by the Russian Foundation of
Basic Research, project 15-07-02984.}}
\author{V. Yu. Korolev\thanks{Faculty of Computational Mathematics
and Cybernetics, Lomonosov Moscow State University; Institute for
Informatics Problems, Informatics and Control Federal Research
Center, vkorolev@cs.msu.su}, A. V. Dorofeeva\thanks{Faculty of
Computational Mathematics and Cybernetics, Lomonosov Moscow State
University, dorofeyeva$\_$sasha@mail.ru}}
\date{}

\begin{document}

\maketitle

{\bf Abstract:} Bounds of the accuracy of the normal approximation
to the distribution of a sum of independent random variables are
improved under relaxed moment conditions, in particular, under the
absence of moments of orders higher than the second. These results
are extended to Poisson-binomial, binomial and Poisson random sums.
Under the same conditions, bounds are obtained for the accuracy of
the approximation of the distributions of mixed Poisson random sums
by the corresponding limit law. In particular, these bounds are
constructed for the accuracy of approximation of the distributions
of geometric, negative binomial and Poisson-inverse gamma (Sichel)
random sums by the Laplace, variance gamma and Student
distributions, respectively. All absolute constants are written out
explicitly.

\smallskip


\normalsize


\section{Introduction}

\subsection{The history of the problem and aims of the paper}

Let $X_1,X_2,\ldots $ be independent random variables with $\E
X_i=0$ and $0<\E X_i^2\equiv\sigma_i^2<\infty$, $i=1,2,\ldots $ For
$n\in\N$ denote $S_n=X_1+\ldots +X_n$, $B_n^2=\sigma_1^2+\ldots
+\sigma_n^2$. Let $\Phi(x)$ be the standard normal distribution
function,
$$
\Phi(x)=\frac{1}{\sqrt{2\pi}}\int\limits_{-\infty}^{x}e^{-z^2/2}dz,\
\ \ \ x\in\R.
$$
Denote
$$
\Delta_n=\sup_x|\P(S_n<xB_n)-\Phi(x)|.
$$
Let $\mathcal{G}$ be the class of real functions $g(x)$ of $x\in\R$
such that\vspace{-2mm}
\begin{itemize}
\item $g(x)$ is even;\vspace{-2mm}
\item $g(x)$ is nonnegative for all $x$ and $g(x)>0$ for $x>0$;\vspace{-2mm}
\item $g(x)$ and $x/g(x)$ do not decrease for $x>0$.\vspace{-1mm}
\end{itemize}

In 1963 M. Katz \cite{Katz1963} proved that, whatever
$g\in\mathcal{G}$ is, if the random variables $X_1,X_2,\ldots $ are
identically distributed and $\E X_1^2g(X_1)<\infty$, then there
exists a finite positive constant $C_1$ such that
$$
\Delta_n\le C_1\cdot\frac{\E X_1^2g(X_1)}{\sigma_1^2
g\big(\sigma_1\sqrt{n}\big)}.\eqno(1)
$$
In 1965 this result was generalized by V.\,V.\,Petrov
\cite{Petrov1965} to the case of non-identically distributed
summands (also see \cite{Petrov1972}): whatever $g\in\mathcal{G}$
is, if $\E X_i^2g(X_i)<\infty$, $i=1,\ldots ,n$, then there exists a
finite positive constant $C_2$ such that
$$
\Delta_n\le \frac{C_2}{B_n^2g(B_n)}\sum_{i=1}^n\E
X_i^2g(X_i).\eqno(2)
$$

Everywhere in what follows the symbol $\I(A)$ will denote the
indicator function of an event $A$. For $\varepsilon\in(0,\infty)$
denote
$$
L_n(\varepsilon)=\frac{1}{B_n^2}\sum_{i=1}^n\E X_i^2\I(|X_i|\ge
\varepsilon B_n),\ \ \
M_n(\varepsilon)=\frac{1}{B_n^3}\sum_{i=1}^n\E |X_i|^3\I(|X_i|<
\varepsilon B_n).
$$
In 1966 L.\,V.~Osipov \cite{Osipov1966} proved that there exists a
finite positive absolute constant $C_3$ such that for any
$\varepsilon\in(0,\infty)$
$$
\Delta_n\le C_3\big[L_n(\varepsilon)+M_n(\varepsilon)\big]\eqno(3)
$$
(also see \cite{Petrov1987}, Chapt V, Sect.\,3, theorem 7). This
inequality is of special importance. Indeed, it is easy to see that
$$
M_n(\varepsilon)\le\frac{\varepsilon}{B_n^2}\sum_{i=1}^{n}\E
X_i^2\I(|X_i|< \varepsilon B_n)\le\varepsilon.
$$
Hence, from (3) it follows that for any $\varepsilon\in(0,\infty)$
$$
\Delta_n\le C_3\big(\varepsilon+L_n(\varepsilon)\big).\eqno(4)
$$
But, as is well known, the Lindeberg condition
$$
\lim_{n\to\infty}L_n(\varepsilon)=0\ \ \text{for any
}\varepsilon\in(0,\infty)
$$
is a {\it criterion} of convergence in the central limit theorem.
Therefore, in terminology proposed by
V.\,M.\,Zolotarev~\cite{Zolotarev1997}, bound (4) is {\it natural},
since it relates the convergence {\it criterion} with the
convergence {\it rate} and its heft-hand and right-hand sides
converge to zero or diverge simultaneously.

In 1968 inequality (3) in a somewhat more general form was re-proved
by W.\,Feller \cite{Feller1968}, who used the method of
characteristic functions to show that $C_3\le6$.

A special case of (3) is the inequality
$$
\Delta_n\le C_3'\big[L_n(1)+M_n(1)].\eqno(5)
$$
In the book \cite{Petrov1972} it was demonstrated that $C_3\le
2C_3'$.

For identically distributed summands inequality (5) takes the form
$$
\Delta_n\le \frac{C_4}{\sigma_1^2}{\sf
E}X_1^2\min\bigg\{1,\,\frac{|X_1|}{\sigma_1\sqrt{n}}\bigg\}.\eqno(6)
$$

In the papers \cite{Paditz1980, Paditz1984} L.~Paditz showed that
the constant $C_4$ can be bounded as $C_4<4.77$. In 1986 in the
paper \cite{Paditz1986} he noted that with the account of lemma 12.2
from \cite{Bhat1982}, using the technique developed in
\cite{Paditz1980, Paditz1984}, the upper bound for $C_4$ can be
lowered to $C_4<3.51$.

In 1984 A.\,Barbour and P.\,Hall \cite{BarbourHall1984} proved
inequality (5) by Stein's method and, citing Feller's result
mentioned above, stated that the method they used gave only the
bound $C_3'\le 18$ (although the paper itself contains only the
proof of the bound $C_3'\le22$). In 2001 L.\,Chen and K.\,Shao
published the paper \cite{ChenShao2001} containing no references to
Paditz' papers \cite{Paditz1980, Paditz1984, Paditz1986} in which
the proved inequality (5) by Stein's method with the absolute
constant $C_3'=4.1$.

In 2011 V.\,Yu.\,Korolev and S.\,V.\,Popov  \cite{KP2011_3} showed
that there exist universal constants $C_1$ and $C_2$ which do not
depend on a particular form of $g\in\mathcal{G}$, such that
inequalities (1), (2), (5) and (6) are valid with $C_1=C_4\le
3.0466$ and $C_2=C_3'\le 3.1905$. This result was later improved by
the same authors in the papers \cite{KorolevPopov2011,
KorolevPopovDAN}, where it was shown that $C_1=C_2=C_4=C_3'\le
2.011$.

Moreover, in the paper \cite{KorolevPopovDAN} lower bounds were
established for the universal constants $C_1$ and $C_2$. Namely, let
$g$ be an arbitrary function from the class $\mathcal{G}$. Denote by
$\mathcal{H}_g$ the set of all random variables $X$ satisfying the
condition ${\sf E}X^2g(X)<\infty$. Denote
$$
C^*=\sup_{g\in\mathcal{G}}\sup_{{X_i\in\mathcal{H}_g,}\atop{i=1,\ldots
,n}}\frac{\Delta_n B_n^2g(B_n)}{\sum_{i=1}^n\E X_i^2g(X_i)}.
$$
It is easily seen that $C^*$ is the least possible value of the
absolute constant $C_2$ that provides the validity of inequality (2)
for all functions $g\in\mathcal{G}$ at once. In the paper
\cite{KorolevPopovDAN} it was proved that
$$
C^*\ge\sup_{z>0}\Big|\frac{1}{1+z^2}-\Phi(-z)\Big|=0.54093\ldots
$$

The aim of the present paper is to improve and extend the results
mentioned above. First, we will show that one can take $C_3=C_3'$.
Second, we will sharpen the upper bounds of the absolute constants
mentioned above. Third, we will extend these results to
Poisson-binomial, binomial and Poisson random sums. Under the same
conditions, bounds will be obtained for the accuracy of the
approximation of the distributions of mixed Poisson random sums by
the corresponding limit law. In particular, we will construct these
bounds for the accuracy of approximation of the distributions of
geometric, negative binomial and Poisson-inverse gamma (Sichel)
random sums by the Laplace, variance gamma and Student
distributions, respectively. All absolute constants will be written
out explicitly.

Along with purely theoretical motivation to sharpen and generalize
known results, there is a somewhat practical interest in the
problems considered below. Poisson-binomial, binomial and mixed
Poisson (first of all, geometric) random sums are widely used as
stopped-random-walk models in many fields such as financial
mathematics (Cox--Ross--Rubinstein binomial random walk model for
option pricing \cite{CoxRossRubinstein1979}), insurance (Poisson
random sums as total claim size in dynamic collective risk models
\cite{Cramer1930}, binomial random sums as total claim size in
static portfolio risk models, geometric sums in the
Pollaczek--Khinchin--Beekman representation of the ruin probability
within the framework of the classical risk process
\cite{Kalashnikov1997}), reliability theory for modeling rare events
\cite{Kalashnikov1997}. It is now a tradition to admit that the
distributions of elementary jumps of these random walks may have
very heavy tails. The problems considered in the present paper
correspond to the situation where the tails may be as heavy as
possible for the normal approximation to be still adequate.
Moreover, the bounds obtained in this paper partly give an answer to
the questions how heavy these tails can be for the normal
approximation (or scale-mixed normal approximation) to be
reasonable.

The paper is organized as follows. In Section 1.1 we prove that in
inequalities (1)--(5) the absolute constants coincide and that the
values of these constants are determined by that of $C_3'$. In
Section 2 the upper bound of $C_3'$ is sharpened. In Section 3 the
analogs of inequalities (1), (2), (3) and (6) are proved for
Poisson-binomial and binomial random sums. In Section 4 the results
obtained in Section 3 are used to construct the analogs of (1) and
(6) for Poisson random sums. The results of Section 4 are used in
Section 5 to obtain bounds for the accuracy of the approximation of
the distributions of mixed Poisson random sums by the corresponding
limit law. In particular, here these bounds are constructed for the
accuracy of approximation of the distributions of geometric,
negative binomial and Poisson-inverse gamma random sums by the
Laplace, variance gamma and Student distributions, respectively.

\subsection{On the coincidence of the absolute constants in inequalities (1)--(5)}

The main result of this section is the following statement.

\smallskip

{\sc Lemma 1.} {\it For any $\varepsilon\in(0,\infty)$}
$$
L_n(1)+M_n(1)\le L_n(\varepsilon)+M_n(\varepsilon).\eqno(7)
$$

\smallskip

{\sc Proof}. For $\varepsilon=1$ the statement is trivial. Let
$\varepsilon<1$. Then
$$
L_n(1)+M_n(1)=L_n(\varepsilon)+M_n(\varepsilon)+
$$
$$
+ \frac{1}{B_n^3}\sum_{j=1}^n|X_j|^3\I(\varepsilon B_n\le|X_j|<
B_n)-\frac{1}{B_n^2}\sum_{j=1}^n{\sf E}X_j^2\I(\varepsilon
B_n\le|X_j|< B_n).
$$
But
$$
\frac{1}{B_n^3}\sum_{j=1}^n|X_j|^3\I(\varepsilon B_n\le|X_j|< B_n)-
\frac{1}{B_n^2}\sum_{j=1}^n{\sf E}X_j^2\I(\varepsilon B_n\le|X_j|<
B_n)\le
$$
$$
\le\frac{1}{B_n^2}\sum_{j=1}^n{\sf E}X_j^2\I(\varepsilon
B_n\le|X_j|< B_n)-\frac{1}{B_n^2}\sum_{j=1}^n{\sf
E}X_j^2\I(\varepsilon B_n\le|X_j|< B_n)=0,
$$
therefore, in the case $\varepsilon<1$ inequality (7) is proved.

Let now $\varepsilon>1$. Then
$$
L_n(1)+M_n(1)=L_n(\varepsilon)+M_n(\varepsilon)+
$$
$$
+ \frac{1}{B_n^2}\sum_{j=1}^n{\sf E}X_j^2\I(B_n\le|X_j|<\varepsilon
B_n)-\frac{1}{B_n^3}\sum_{j=1}^n|X_j|^3\I(B_n\le|X_j|<\varepsilon
B_n).
$$
But
$$
\frac{1}{B_n^2}\sum_{j=1}^n{\sf E}X_j^2\I(B_n\le|X_j|<\varepsilon
B_n)-\frac{1}{B_n^3}\sum_{j=1}^n|X_j|^3\I(B_n\le|X_j|<\varepsilon
B_n)\le
$$
$$
\le\frac{1}{B_n^2}\sum_{j=1}^n{\sf E}X_j^2\I(B_n\le|X_j|<\varepsilon
B_n)-\frac{1}{B_n^2}\sum_{j=1}^n{\sf
E}X_j^2\I(B_n\le|X_j|<\varepsilon B_n)=0,
$$
that is, the statement of the lemma holds for $\varepsilon>1$. as
well.

\smallskip

{\sc Corollary 1.} {\it The absolute constants in inequalities
$(3)$, $(4)$, $(5)$ and $(6)$ can be taken identical, that is, if
inequality $(5)$ holds with $C_3'\le C_0$, then inequalities $(3)$,
$(4)$ and $(6)$ hold with $C_3\le C_0$ и $C_4\le C_0$.}

\smallskip

{\sc Remark 1.} In the paper \cite{KorolevPopovDAN} it was shown
that if inequality $(5)$ holds with $C_3'\le C_0$, then inequalities
$(1)$ and $(2)$ hold with $C_i\le C_0$, $i=1,2$.

\smallskip

So, in the evaluation of the constants in the above inequalities,
the constant $C_3'$ in inequality (5) plays the determining role: if
a particular upper bound $C_3'\le C_0$ is known, then in all the
rest inequalities (1)--(4) and (6) one can let $C_i\le C_0$,
$i=1,2,3,4$. That is the reason for us to focus on sharpening the
upper bound for $C_3'$.

\section{Sharpening of the upper bound for the constant $C_3'$}

\subsection{Auxiliary results}

For $x\ge 0$, $n\in\N$ and $i=1,\ldots,n$ denote
$$
Y_i(x)=B_n^{-1}X_i \I \big(|X_i|<(1+x)B_n\big),\ \  \ Y_i=Y_i(0),\ \
\ W_n(x)=\sum_{i=1}^nY_i(x),\ \ \ \ W_n=W_n(0).
$$
Since $\E X_i=0$, we have
$$
|\E X_i\I \big(|X_i|<(1+x)B_n\big)|=|\E X_i\I \big(|X_i|\ge
(1+x)B_n\big)|.\eqno(8)
$$
By the definition of the random variables $Y_i(x)$ the relation
$$
\sum_{i=1}^n\E Y_i^2(x)\le \frac{1}{B_n^2}\sum_{i=1}^n\E X_i^2=
1\eqno(9)
$$
holds. Denote
$$
K=\frac{17+7\sqrt{7}}{27}<1.3156.
$$

\smallskip

{\sc Lemma 2}. {\it $1^{\circ}$. For any $n\in\N$, $x\ge0$ and
$p\in[1,\,K]$ there holds the inequality
$$
\sum_{i=1}^n\E |Y_i(x)-\E Y_i(x)|^3\le \min\Big\{K
M_n(1+x),\,pM_n(1+x)+\frac{(5-p)\,L_n(1+x)}{1+x}\Big\}.
$$

\noindent $2^{\circ}$. For any $n\in\N$ and $x\ge0$ there hold the
inequalities
$$
1-2L_n(1+x)\le\D W_n(x)\le 1.
$$

\noindent $3^{\circ}$. Let $M_n(1)=\gamma L_n(1)$, $\gamma\ge 0$.
Then for any $n\in\N$ there holds the inequality}
$$
\sum_{i=1}^n\E |Y_i-\E Y_i|^3\le L_n(1)
\min\left\{K\gamma,\,\gamma+4\right\}.
$$

\smallskip

The {\sc proof} based on the results of \cite{Hoeffding1948,
Zolotarev1986} and \cite{NefedovaShevtsova2011} was given in
\cite{KorolevPopov2011}.

\smallskip

{\sc Lemma 3}. {\it $1^{\circ}$. Let $q>0$. Then
$$
\sup_x|\Phi(qx)-\Phi(x)|=\bigg|\Phi\bigg(q\sqrt{\frac{\ln
q^2}{q^2-1}}\bigg)-\Phi\bigg(\sqrt{\frac{\ln
q^2}{q^2-1}}\bigg)\bigg|\le
$$
$$
\le\sqrt{\frac{(q-1)\ln q}{\pi(q+1)}}\exp\Big\{-\min(1,\,q)\frac{\ln
q}{q^2-1}\Big\} \le\frac{1}{\sqrt{2\pi
e}}\Big(\max\Big\{q,\,\frac{1}{q}\Big\}-1\Big).
$$

\noindent $2^{\circ}$. Let $a\in\R$. Then}
$$
\sup_x|\Phi(x+a)-\Phi(x)|=2\Phi\Big(\frac{\abs
a\abs}{2}\Big)-1\le\frac{\abs a\abs}{\sqrt{2\pi}}.
$$

\smallskip

The elementary {\sc proof} of this lemma is based on the Lagrange
formula and the easily verifiable fact: if $F(x)$ and $G(x)$ are two
differentiable distribution functions, then $\sup\limits_x|F(x)-G(x)|$ is
attained at those points $x$, where $F'(x)=G'(x)$ (also see
\cite{Petrov1972}, p. 143).

\smallskip

{\sc Lemma 4}. {\it Assume that $L_n(1)\le A$ for some
$A\in(0,\,\frac12)$. Let
$$
B(A)=\frac{2}{(1+\sqrt{1-2A})\sqrt{1-2A}}.
$$
Then}
$$
1\le\frac{1}{\sqrt{\D W_n}}\le 1+B(A)L_n(1).
$$

\smallskip

For the {\sc proof} see \cite{KorolevPopov2011}.

\smallskip

{\sc Lemma 5}. {\it Let $X$ be a random variable with $\E
X^2<\infty$. Then}
$$
\sup_x\bigg|\P\bigg(\frac{X-\E X}{\sqrt{\D
X}}<x\bigg)-\Phi(x)\bigg|\le
\sup_{z>0}\Big|\frac{1}{1+z^2}-\Phi(-z)\Big|=0.54093\ldots
$$

\smallskip

For the {\sc proof} see, e. g., the book \cite{Bhat1982} and the
papers \cite{Kondrik, KorolevPopovDAN}.

\subsection{General case}

{\sc Theorem 1}. {\it Let $n\in\mathbb{N}$, the random variables
$X_1,\ldots,X_n$ be independent, $\E X_i=0$ and $0<\E X_i^2<\infty$,
$i=1,\ldots,n$. Let $\gamma=M_n(1)/L_n(1)$. Then there exists a
finite positive number $C_1(\gamma)$ depending only on $\gamma$ such
that
$$
\Delta_n\le (1+\gamma)C_1(\gamma)L_n(1).
$$
Moreover, the upper bounds for $C_1(\gamma)$ are presented in table
$1$.}

\renewcommand{\tablename}{Table}

\renewcommand{\baselinestretch}{1.2}
\begin{table}[h]
\centering \small
  \begin{tabular}{||c|c||c|c||c|c||}
  \hline
  $\gamma$ & $C_1(\gamma)\le$ & $\gamma$ & $C_1(\gamma)\le$ & $\gamma$ & $C_1(\gamma)\le$ \\
\hline
$\gamma \ge 0$           & 1.8627 & $\gamma \ge 1$   & 1.5605 & $\gamma \ge 10$  & 0.9393 \\
$\gamma \ge 0.1$         & 1.8587 & $\gamma \ge 2$   & 1.3488 & $\gamma \ge 100$ & 0.6067 \\
$\gamma \ge 0.5$         & 1.7244 & $\gamma \ge 5$   & 1.0836 & $\gamma\rightarrow\infty$ & 0.5583 \\
\hline
\end{tabular}
\caption{Upper bounds for $C_1(\gamma)$.}
\end{table}
\renewcommand{\baselinestretch}{1.0}

\normalsize

\smallskip

{\sc Corollary 2}. {\it Under the conditions of theorem $1$,
inequalities $(2)-(5)$ hold with $C_2=C_3=C_3'\le 1.8627$.}

\smallskip

{\sc Proof} of theorem 1. For any $y\in\R$ the event $\{S_n<yB_n\}$
implies the event
$$
\{W_n<y\}\cup\{|X_1|\ge B_n\}\cup\ldots\cup\{|X_n|\ge B_n\},
$$
whereas the event $\{W_n<y\}$ implies the event
$$
\{S_n<yB_n\}\cup\{|X_1|\ge B_n\}\cup\ldots\cup\{|X_n|\ge B_n\}.
$$
Therefore,
$$
\sup_y|\P(S_n<yB_n)-\P(W_n<y)|\le\sum_{i=1}^n\P(|X_i|\ge B_n).
$$
Hence, for any $y\in\R$
$$
\Delta_n\le Q_1+Q_2+Q_3,\eqno(10)
$$
where
$$
Q_1=\sup_y\left|\P\left(\frac{W_n-\E W_n}{\sqrt{\D W_n}}<\frac{y-\E
W_n}{\sqrt{\D W_n}}\right)-\Phi\left(\frac{y-\E W_n}{\sqrt{\D
W_n}}\right)\right|,
$$
$$
Q_2=\sup_y\left|\Phi\left(\frac{y-\E W_n}{\sqrt{\D
W_n}}\right)-\Phi(y)\right|,\ \ \ \ Q_3=\sum_{i=1}^n\P(|X_i|\ge
B_n).
$$
Consider $Q_1$. By virtue of the Berry--Esseen inequality with the
best known upper bound of the absolute constant (see
\cite{Shevtsova2013}) we have
$$
Q_1\le\frac{0.5583}{\big(\D W_n\big)^{3/2}}\sum_{i=1}^n\E|Y_i-\E
Y_i|^3.
$$
Assume that $L_n(1)\le A<\frac12$. Then in accordance with
statements $2^{\circ}$ and $3^{\circ}$ of lemma 2
$$
Q_1\le\frac{0.5583\cdot
 \min\{K\gamma,\,\gamma+4\}\,L_n(1)}{(1-2A)^{3/2}}.\eqno(11)
$$
Consider $Q_2$. We obviously have
$$
Q_2=\sup_y\left|\Phi\left(\frac{y-\E W_n}{\sqrt{\D
W_n}}\right)-\Phi\left(y-\E W_n\right)+ \Phi\left(y-\E
W_n\right)-\Phi(y)\right|\le
$$
$$
\le\sup_y\left|\Phi\left(\frac{y-\E W_n}{\sqrt{\D
W_n}}\right)-\Phi\left(y-\E W_n\right)\right|+
\sup_y\left|\,\Phi\!\left(y-\E W_n\right)-\Phi(y)\right|=
$$
$$
=\sup_y\left|\Phi\left(\frac{y}{\sqrt{\D
W_n}}\right)-\Phi\left(y\right)\right|+
\sup_y\left|\,\Phi\!\left(y-\E W_n\right)-\Phi(y)\right|\equiv
Q_{21}+Q_{22}.
$$
According to statement $2^{\circ}$ of lemma 2, $\D W_n\le 1$.
Therefore, by virtue of statement $1^{\circ}$ of lemma 3 and lemma
4, there holds the inequality
$$
Q_{21}\le\frac{1}{\sqrt{2\pi e}}\left(\frac{1}{\sqrt{\D
W_n}}-1\right)\le \frac{2L_n(1)}{\sqrt{2\pi
e(1-2A)}(1+\sqrt{1-2A})}.\eqno(12)
$$
Consider $Q_{22}$. By virtue of (8) we have
$$
|\E W_n|=\Big|\sum_{i=1}^n\E Y_i\Big|\le\frac{1}{B_n}\sum_{i=1}^n|\E
X_i\I(|X_i|<B_n)|=\frac{1}{B_n}\sum_{i=1}^n|\E X_i\I(|X_i|\ge
B_n)|\le
$$
$$
\le\frac{1}{B_n}\sum_{i=1}^n\E |X_i|\I(|X_i|\ge
B_n)\le\frac{1}{B_n^2}\sum_{i=1}^n\E X_i^2\I(|X_i|\ge B_n)=L_n(1).
$$
Therefore, by statements $2^{\circ}$ of lemma 2 and $2^{\circ}$ of
lemma 3,
$$
Q_{22}\le\frac{L_n(1)}{\sqrt{2\pi}}.\eqno(13)
$$
From (12) and (13) it follows that
$$
Q_2\le\frac{L_n(1)}{\sqrt{2\pi}}\left(1+\frac{2}{\sqrt{e(1-2A)}(1+\sqrt{1-2A})}\right).\eqno(14)
$$
Finally, by the Markov inequality
$$
Q_3=\sum_{i=1}^n\P(|X_i|\ge B_n)\le\frac{1}{B_n^2}\sum_{i=1}^n\E
X_i^2\I(|X_i|\ge B_n)=L_n(1).\eqno(15)
$$
So, from (10), (11), (14) and (15) we obtain
$$
\Delta_n\le
L_n(1)\left[1\!+\!\frac{1}{\sqrt{2\pi}}\left(1\!+\!\frac{2}{\sqrt{e(1-2A)}(1+\sqrt{1-2A})}\right)\!+\!
\frac{0.5583\cdot\min\{K\gamma,\,\gamma+4\}}{(1-2A)^{3/2}}\right].\eqno(16)
$$
Introduce the function
$$
H_1(\gamma,
A)=1\!+\!\frac{1}{\sqrt{2\pi}}\left(1\!+\!\frac{2}{\sqrt{e(1-2A)}(1+\sqrt{1-2A})}\right)\!+
\!\frac{0.5583\cdot\min\{K\gamma,\,\gamma+4\}}{(1-2A)^{3/2}}.
\eqno(17)
$$
For any $0\le A<\frac{1}{2}$ we have the inequality
$$
\Delta_n\le
L_n(1)\cdot\max\left\{H_1(\gamma,A),\,\frac{0.541}{A}\right\}.
$$
This follows from (16) if $L_n(1)\le A$ and from lemma 5 otherwise.

Now, with the account of the equality
$$
L_n(1)=\frac{(L_n(1)+\gamma
L_n(1))}{1+\gamma}=\frac{L_n(1)+M_n(1)}{1+\gamma},
$$
we have
$$
C_1(\gamma)\le \min_{0\le
A<\frac{1}{2}}\max\left\{\frac{H_1(\gamma,A)}{1+\gamma},\,\frac{0.541}{A(1+\gamma)}\right\}.
$$
The computation by this formula yield the values presented in table
1. Note that the first function of $A$ inside the minimax is
increasing whereas the second one is decreasing. Hence, the value of
the minimax is delivered by the unique solution of the equation
$$
\frac{H_1(\gamma,A)}{1+\gamma}\,=\,\frac{0.541}{A(1+\gamma)}.
$$
For $\gamma>13$ (we have $\gamma+4<K\gamma$) both functions decrease
in $\gamma$, that is, the minimax value decreases. Therefore, the
corresponding part of table 1 is obtained by the evaluation of the
bound for $C_1(\gamma)$ at one point. The part of table 1
corresponding to $0\le\gamma\le 13$ is obtained by numerical
optimization of a finite interval. The theorem is proved.

\subsection{Special cases}

Using the best current upper bound $C_0\le0.4690$ for the absolute
constant in the Berry--Esseen inequality for identically distributed
summands (see \cite{Shevtsova2013}), the following statement can be
obtained in the way similar to the proof of theorem 1.

\smallskip

{\sc Theorem 2}. {\it In addition to the assumptions of theorem $1$,
let the random variables $X_1,X_2,\ldots$ be identically
distributed. Then there exists a finite positive number
$C_2(\gamma)$ depending only on $\gamma$ such
that
$$
\Delta_n\le(1+\gamma)C_2(\gamma)L_n(1).
$$
Moreover, the upper bounds for $C_2(\gamma)$ are presented in table
$2$.}
\renewcommand{\baselinestretch}{1.2}
\begin{table}[h]
\centering \small
  \begin{tabular}{||c|c||c|c||c|c||}
  \hline
  $\gamma$ & $C_2(\gamma)\le$ & $\gamma$ & $C_2(\gamma)\le$ & $\gamma$ & $C_2(\gamma)\le$  \\
\hline
$\gamma \ge 0$           & 1.8546 & $\gamma \ge 1$           & 1.4793 & $\gamma \ge 10$           & 0.8292 \\
$\gamma \ge 0.1$         & 1.8338 & $\gamma \ge 2$           & 1.2540 & $\gamma \ge 100$          & 0.5147 \\
$\gamma \ge 0.5$         & 1.6608 & $\gamma \ge 5$           & 0.9781 & $\gamma\rightarrow\infty$ & 0.4690 \\
\hline
\end{tabular}
\caption{Upper bounds for $C_2(\gamma)$.}
\end{table}
\renewcommand{\baselinestretch}{1.0}

\normalsize

\smallskip

{\sc Proof}. Using the reasoning similar to that used to prove
theorem 1, it is easy to see that
$$
C_2(\gamma)\le \min_{0\le
A<\frac{1}{2}}\max\left\{\frac{H_2(\gamma,A)}{1+\gamma},\,\frac{0.541}{A(1+\gamma)}\right\},\vspace{-2mm}
$$
where
$$
\vspace{-2mm}H_2(\gamma,
A)=1\!+\!\frac{1}{\sqrt{2\pi}}\left(1\!+\!\frac{2}{\sqrt{e(1-2A)}(1+\sqrt{1-2A})}\right)\!+
\!\frac{0.4690\cdot\min\{K\gamma,\,\gamma+4\}}{(1-2A)^{3/2}}.
$$
The computations by these formula yield the values of the upper
bounds for $C_2(\gamma)$ presented in table 2. The theorem is
proved.

\smallskip

{\sc Corollary 3}. {\it Under conditions of theorem $2$,
inequalities $(1)$ and $(6)$ hold with $C_1=C_4\le 1.8546$.}

\smallskip

{\sc Theorem 3}. {\it In addition to the conditions of theorem $1$,
let the random variables $X_1,X_2,\ldots$ have symmetric
distributions. Then there exists a finite positive number
$C_3(\gamma)$ depending only on $\gamma$ such that
$$
\Delta_n\le(1+\gamma)C_3(\gamma)L_n(1).
$$
Moreover, the upper bounds for $C_3(\gamma)$ are presented in table
$3$.}

\renewcommand{\baselinestretch}{1.2}
\begin{table}[h]
\centering \small
  \begin{tabular}{||c|c||c|c||c|c||}
  \hline
  $\gamma$ & $C_3(\gamma)\le$ & $\gamma$ & $C_3(\gamma)\le$ & $\gamma$ & $C_3(\gamma)\le$ \\
\hline
$\gamma \ge 0$           & 1.5769 & $\gamma \ge 1$           & 1.3033 & $\gamma \ge 10$  & 0.7433 \\
$\gamma \ge 0.1$         & 1.5749 & $\gamma \ge 2$           & 1.1115 & $\gamma \ge 100$ & 0.5808 \\
$\gamma \ge 0.5$         & 1.4532 & $\gamma \ge 5$           & 0.8729 & $\gamma\rightarrow\infty$ & 0.5583 \\
\hline
\end{tabular}
\caption{Upper bounds of $C_3(\gamma)$.}
\end{table}
\renewcommand{\baselinestretch}{1.0}

\normalsize

\smallskip

{\sc Corollary 4}. {\it Under the conditions of theorem $3$,
inequalities $(2)-(5)$ hold with $C_2=C_3=C_3'\le 1.5769$.}

\smallskip

The {\sc proof} of theorem 3. In the case under consideration
instead (8) we have
$$
Q_1\le\frac{0.5583M_n(1)}{(1-2A)^{3/2}},
$$
and $Q_{22}=0$, since $\E W_n=0$. Therefore, the bound
$$
\Delta_n\le L_n(1)\bigg(1+\frac{2}{\sqrt{2\pi
e(1-2A)}(1+\sqrt{1-2A})}\bigg)+\frac{0.5583M_n(1)}{(1-2A)^{3/2}}
$$
holds. Thus,
$$
C_3(\gamma)\le \min_{0\le
A<\frac{1}{2}}\max\left\{\frac{H_3(\gamma,A)}{1+\gamma},\,\frac{0.541}{A(1+\gamma)}\right\},
$$
where
$$
H_3(\gamma,A)=1+\frac{2}{\sqrt{2\pi
e(1-2A)}(1+\sqrt{1-2A})}+\frac{0.5583\gamma}{(1-2A)^{3/2}}.
$$
The computations by the above formulas yield the values of the upper
bounds for $C_3(\gamma)$ presented in table 3. The theorem is
proved.

\smallskip

{\sc Theorem 4}. {\it In addition to the conditions of $3$, let the
random variables $X_1,X_2,\ldots$ be identically distributed. Then
there exists a finite positive number $C_4(\gamma)$ depending only
on $\gamma$ such that
$$
\Delta_n\le(1+\gamma)C_4(\gamma)L_n(1).
$$
Moreover, the upper bounds for $C_4(\gamma)$ are presented in table
$4$.}

\renewcommand{\baselinestretch}{1.2}
\begin{table}[h]
\centering \small
  \begin{tabular}{||c|c||c|c||c|c||}
  \hline
  $\gamma$ & $C_4(\gamma)\le$ & $\gamma$ & $C_4(\gamma)\le$ & $\gamma$ & $C_4(\gamma)\le$ \\
\hline
$\gamma \ge 0$           & 1.5645 & $\gamma \ge 1$           & 1.2388 & $\gamma \ge 10$  & 0.6591 \\
$\gamma \ge 0.1$         & 1.5534 & $\gamma \ge 2$           & 1.0373 & $\gamma \ge 100$ & 0.4923 \\
$\gamma \ge 0.5$         & 1.4018 & $\gamma \ge 5$           & 0.7915 & $\gamma\rightarrow\infty$ & 0.4690 \\
\hline
\end{tabular}
\caption{Upper bound for $C_4(\gamma)$.}
\end{table}
\renewcommand{\baselinestretch}{1.0}

\normalsize

\smallskip

{\sc Corollary 5}. {\it Under the conditions of theorem $4$,
inequalities $(1)$ and $(6)$ hold with $C_1=C_4\le 1.5645$.}

\smallskip

{\sc Proof} of theorem 4. In the case under consideration
$$
C_4(\gamma)\le \min_{0\le
A<\frac{1}{2}}\max\left\{\frac{H_4(\gamma,A)}{1+\gamma},\,\frac{0.541}{A(1+\gamma)}\right\},
$$
where
$$
H_4(\gamma,A)=1+\frac{2}{\sqrt{2\pi
e(1-2A)}(1+\sqrt{1-2A})}+\frac{0.4690\gamma}{(1-2A)^{3/2}}.
$$
The computations by the above formulas yield the values of the upper
bounds for $C_4(\gamma)$ presented in table 4. The theorem is
proved.

\section{The accuracy of the normal approximation to the distributions of Poisson-binomial random sums}

From this point on let $X_1,X_2,\ldots$ be independent {\it
identically distributed} random variables with $\E X_i=0$ and $0<\E
X_i^2\equiv \sigma^2<\infty$. Let $p_j\in(0,1]$ be arbitrary
numbers, $j=1,2,\ldots$. For $n\in\N$ denote
$\theta_n=p_1+\ldots+p_n$, $\vp_n=(p_1,\ldots,p_n)$. The
distribution of the random variable
$$
N_{n,\vp_n}=\xi_1+\ldots+\xi_n,
$$
where $\xi_1,\ldots,\xi_n$ are independent random variables such
that
$$
\xi_j=\begin{cases}1 & \text{ with probability } p_j,\cr
                    0 & \text{ with probability } 1-p_j,
      \end{cases},\ \ \ j=1,\ldots,n,
$$
is usually called Poisson-binomial distribution with parameters
$n;\vp_n$. Assume that for each $n\in\N$ the random variables
$N_{n,\vp_n},X_1,X_2,\ldots$ are jointly independent. The main
objects considered in this section are {\it Poisson-binomial random
sums} of the form
$$
S_{N_{n,\vp_n}}=X_1+\ldots+X_{N_{n,\vp_n}}.
$$
As this is so, if $N_{n,\vp_n}=0$, then we assume
$S_{N_{n,\vp_n}}=0$.

For $j\in\mathbb{N}$ introduce the random variables $\wX_j$ by
setting
$$
\wX_j=\begin{cases}X_j & \text{ with probability } p_j,\cr 0 &
\text{ with probability } 1-p_j.\end{cases}
$$
If the common distribution function of the random variables $X_j$ is
denoted $F(x)$ and the distribution function with a single unit jump
at zero is denoted $E_0(x)$, then, as is easily seen,
$$
{\sf P}\big(\wX_j<x\big)=p_jF(x)+(1-p_j)E_0(x),\ \ \ \
x\in\mathbb{R},\ j\in\mathbb{N}.
$$
It is obvious that ${\sf E}\wX_j=0$,
$$
{\sf D}\wX_j={\sf E}\wX_j^2=p_j\sigma^2.\eqno(18)
$$
In what follows the symbol $\eqd$ will denote coincidence of
distributions.

\smallskip

{\sc Lemma 6}. {\it For any $n\in\N$ and $p_j\in(0,1]$}
$$
S_{N_{n,\vp_n}}\eqd \wX_1+\ldots+\wX_n, \eqno(19)
$$
{\it where the random variables on the right-hand side of $(19)$ are
independent.}

\smallskip

{\sc Proof.} The characteristic
functions of the left-hand and right-hand sides of (19) have the following forms
$$\varphi_{S_{N_{n,\vp_n}}}(t) = \sum\limits_{k = 0}^n\varphi_{X_1+ \ldots + X_k}(t)\P(N_{n,\vp_n} = k)\;\;\;\;\text{and}\;\;\;\;
\varphi_{\tilde X_1+ \ldots +\tilde X_n}(t) = \prod\limits_{j = 1}^n[p_j\varphi_{X_j}(t) + (1 - p_j)].$$

It suffices to make sure that the characteristic
functions of the left-hand and right-hand sides of (19) coincide.

We will use the method of mathematical induction. Basis: $n  =
1$. 
$$p_1\varphi_{X_1}(t) + (1 - p_1) = p_1\varphi_{X_1}(t) + (1 - p_1).$$

Inductive step: we show that if the characteristic
functions of the left-hand and right-hand sides of (19) coincide with $n = m$, then they also coincide with $n = m + 1$. 
$$\prod\limits_{j = 1}^{m + 1}[p_j\varphi_{X_j}(t) +(1 - p_j)] = (\prod\limits_{j = 1}^{m}[p_j\varphi_{X_j}(t) + (1 - p_j)]) (p_{m + 1}\varphi_{X_{m + 1}}(t) + (1 - p_{m + 1})) =$$
$$= (1 - p_{m + 1})\sum\limits_{k = 0}^m\varphi_{X_1+ \ldots + X_k}(t)\P(N_{m,\vp_m} = k) + p_{m + 1}\varphi_{X_{m + 1}}(t)\sum\limits_{k = 0}^m\varphi_{X_1+ \ldots + X_k}(t)\P(N_{m,\vp_m} = k) = $$
$$= (1 - p_{m + 1})\sum\limits_{k = 0}^m\varphi_{X_1+ \ldots + X_k}(t)\P(N_{m,\vp_m} = k) + p_1\ldots p_{m + 1}\varphi_{X_1 +\ldots+ X_{m + 1}} + $$
$$+ p_{m + 1}\varphi_{X_m + 1}(t)\sum\limits_{k = 0}^{m - 1}\varphi_{X_1+ \ldots + X_k}(t)\P(N_{m,\vp_m} = k) =$$
$$= (\text{note that } \varphi_{X_1}(t) =\ldots =\varphi_{X_{m + 1}}(t) \text{ and transform the last term}) =$$
$$ = (1 - p_{m + 1})\sum\limits_{k = 0}^m\varphi_{X_1+ \ldots + X_k}(t)\P(N_{m,\vp_m} = k) + p_1\ldots p_{m + 1}\varphi_{X_1 +\ldots+ X_{m + 1}} + $$
$$+ p_{m + 1}\sum\limits_{k = 0}^{m}\varphi_{X_1+ \ldots + X_k}(t)\P(N_{m,\vp_m} = k - 1).$$

On the other hand,
$$\sum\limits_{k = 0}^{m + 1}\varphi_{X_1+ \ldots + X_k}(t)\P(N_{m+1,\vp_{m + 1}} = k) = $$
$$= \sum\limits_{k = 0}^{m}\varphi_{X_1+ \ldots + X_k}(t)\P(N_{m+1,\vp_{m + 1}} = k) + p_1\ldots p_{m + 1}\varphi_{X_1 +\ldots+ X_{m + 1}} = $$
$$= \sum\limits_{k = 0}^{m}\varphi_{X_1+ \ldots + X_k}(t)\P(\lbrace N_{m,\vp_{m}} = k \cap \xi_{m + 1} = 0\rbrace \cup \lbrace N_{m,\vp_{m}} = k- 1 \cap \xi_{m + 1} = 1\rbrace) + $$
$$+ p_1\ldots p_{m + 1}\varphi_{X_1 +\ldots+ X_{m + 1}} = $$
$$= (1 - p_{m + 1})\sum\limits_{k = 0}^{m}\varphi_{X_1+ \ldots + X_k}(t)\P( N_{m,\vp_{m}} = k) + p_{m + 1}\sum\limits_{k = 0}^{m}\varphi_{X_1+ \ldots + X_k}(t)\P( N_{m,\vp_{m}} = k - 1) +$$
$$+ p_1\ldots p_{m + 1}\varphi_{X_1 +\ldots+ X_{m + 1}}.$$
Note that the right-hand sides of the above chain of equalities coincide.
The lemma is proved.

\smallskip

With the account of (18) and (19) it is easy to notice that
$$
{\sf D}S_{N_{n,\vp_n}}=\theta_n\sigma^2.\eqno(20)
$$
Denote
$$
\Delta_{n,\vp_n}=\sup_x\big|{\sf
P}\big(S_{N_{n,\vp_n}}<x\sigma\sqrt{\theta_n}\big)-\Phi(x)\big|.
$$

\smallskip

{\sc Theorem 5}. {\it For any $n\in\N$ and $p_j\in(0,1]$, $j\in\N$,}
$$
\Delta_{n,\vp_n}\le\frac{1.8627}{\sigma^2}\,{\sf
E}X_1^2\min\bigg\{1,\,\frac{|X_1|}{\sigma\sqrt{\theta_n}}\bigg\}.
$$

\smallskip

{\sc Proof}. By virtue of lemma 6 and relation (20) we have
$$
\Delta_{n,\vp_n}=\sup_x\big|{\sf P}\big(\wX_1+\ldots+\widetilde
X_n<x\sigma\sqrt{\theta_n}\big)-\Phi(x)\big|,
$$
and for the latter expression we can use the bound given in theorem
1:
$$
\sup_x\big|{\sf
P}\big(\wX_1+\ldots+\wX_n<x\sigma\sqrt{\theta_n}\big)-\Phi(x)\big|\le
$$
$$
\le1.8627\bigg[\frac{1}{\sigma^2\theta_n}\sum_{j=1}^n{\sf
E}\wX_j^2\I\big(|\wX_j|\ge\sigma\sqrt{\theta_n}\big)+\frac{1}{\sigma^3\theta_n^{3/2}}\sum_{j=1}^n{\sf
E}|\wX_j|^3\I\big(|\wX_j|<\sigma\sqrt{\theta_n}\big)\bigg]=
$$
$$
=1.8627\bigg[\frac{1}{\sigma^2\theta_n}\sum_{j=1}^np_j{\sf
E}X_j^2\I\big(|X_j|\ge\sigma\sqrt{\theta_n}\big)+\frac{1}{\sigma^3\theta_n^{3/2}}\sum_{j=1}^np_j{\sf
E}|X_j|^3\I\big(|X_j|<\sigma\sqrt{\theta_n}\big)\bigg]=
$$
$$
=1.8627\bigg[\frac{{\sf
E}X_1^2\I\big(|X_1|\ge\sigma\sqrt{\theta_n}\big)}{\sigma^2\theta_n}\sum_{j=1}^np_j+\frac{{\sf
E}|X_1|^3\I\big(|X_1|<\sigma\sqrt{\theta_n}\big)}{\sigma^3\theta_n^{3/2}}\sum_{j=1}^np_j\bigg]=
$$
$$
=1.8627\bigg[\frac{1}{\sigma^2}{\sf
E}X_{1}^2\I\big(|X_{1}|\ge\sigma\sqrt{\theta_n}\big)+\frac{1}{\sigma^3\sqrt{\theta_n}}{\sf
E}|X_{1}|^3\I\big(|X_{1}|<\sigma\sqrt{\theta_n}\big)\bigg]=
$$
$$
= \frac{1.8627}{\sigma^3\sqrt{\theta_n}}\,{\sf
E}X_1^2\min\big\{\sigma\sqrt{\theta_n},\,|X_1|\big\}=\frac{1.8627}{\sigma^2}\,{\sf
E}X_1^2\min\bigg\{1,\,\frac{|X_1|}{\sigma\sqrt{\theta_n}}\bigg\},
$$
Q. E. D.

\smallskip

{\sc Theorem 6.} {\it Under the conditions of theorem $5$, whatever
function $g\in\mathcal{G}$ is such that $\E X_1^2g(X_1)<\infty$,
there holds the inequality}
$$
\Delta_{n,\vp_n}\le 1.8627\frac{{\sf
E}X_1^2g(X_1)}{\sigma^2g(\sigma\sqrt{\theta_n})}.
$$

\smallskip

{\sc Proof.} Let $g$ be an arbitrary function from the class
$\mathcal{G}$. With the account of the properties of a function
$g\in\mathcal{G}$ it is easy to see that
$$
\E X_1^2\I(|X_1|\ge\sigma\sqrt{\theta_n})=\E
X_1^2\frac{g(X_1)}{g(X_1)}\I(|X_1|\ge \sigma\sqrt{\theta_n})\le
\frac{1}{g(\sigma\sqrt{\theta_n})}\E
X_1^2g(X_1)\I(|X_1|\ge\sigma\sqrt{\theta_n})\eqno(21)
$$
and
$$
\E X_1^3\I(|X_1|<\sigma\sqrt{\theta_n})=\E
X_1^2g(X_1)\frac{|X_1|}{g(X_1)}\I(|X_1|< \sigma\sqrt{\theta_n})\le
\frac{\sigma\sqrt{\theta_n}}{g(\sigma\sqrt{\theta_n})}\E
X_1^2g(X_1)\I(|X_1|<
\sigma\sqrt{\theta_n}).\eqno(22)
$$
Substituting these estimates into the inequality
$$
\Delta_{n,\vp_n}\le 1.8627\bigg[\frac{1}{\sigma^2}{\sf
E}X_{1}^2\I\big(|X_{1}|\ge\sigma\sqrt{\theta_n}\big)+\frac{1}{\sigma^3\sqrt{\theta_n}}{\sf
E}|X_{1}|^3\I\big(|X_{1}|<\sigma\sqrt{\theta_n}\big)\bigg]\eqno(23)
$$
obtained in the proof of theorem 5, we have
$$
\Delta_{n,\vp_n}\le
\frac{1.8627}{\sigma^2g(\sigma\sqrt{\theta_n})}\big[\E
X_1^2g(X_1)\I(|X_1|\ge\sigma\sqrt{\theta_n})+ \E
X_1^2g(X_1)\I(|X_1|< \sigma\sqrt{\theta_n})\big]=1.8627\frac{{\sf
E}X_1^2g(X_1)}{\sigma^2g(\sigma\sqrt{\theta_n})}.
$$
The theorem is proved.

\smallskip

In particular, if $p_1=p_2=\ldots=p$, then the Poisson-binomial
distribution with parameters $n\in\N$ and $\vp_n$ becomes the
classical binomial distribution with parameters $n$ and $p$:
$$
N_{n,\vp_n}\eqd N_{n,p},\ \ \ {\sf
P}(N_{n,p}=k)=C_n^kp^k(1-p)^{n-k},\ \ \ k=0,\ldots,n.
$$
In this case $\theta_n=np$, so that ${\sf D}S_{N_{n,p}}=np\sigma^2$.
Denote
$$
\Delta_{n,p}=\sup_x\big|{\sf
P}\big(S_{N_{n,p}}<x\sigma\sqrt{np}\big)-\Phi(x)\big|.
$$
Estimates of the accuracy of the normal approximation to the
distributions of binomial random sums (under traditional conditions
of the existence of the third moments of summands) were considered
in the paper \cite{Sunklodas}, where a conventional approach was
used which is based on the direct application of the total
probability formula and does not involve representation (19). Hence,
in \cite{Sunklodas} estimates were obtained with the structure far
from being optimal, containing unnecessary terms and unreasonably
large values of absolute constants.

Theorems 2 and 5 imply

\smallskip

{\sc Corollary 6}. {\it For any $n\in\N$ and $p\in(0,1]$}
$$
\Delta_{n,p}\le\frac{1.8546}{\sigma^2}\,{\sf
E}X_1^2\min\bigg\{1,\,\frac{|X_1|}{\sigma\sqrt{np}}\bigg\}.
$$

\smallskip

Theorems 2 and 6 imply

\smallskip

{\sc Corollary 7}. {\it Under the conditions of theorem $5$,
whatever function $g\in\mathcal{G}$ is such that $\E
X_1^2g(X_1)<\infty$, for any $n\in\N$ and $p\in(0,1]$ there holds
the inequality}
$$
\Delta_{n,\vp_n}\le 1.8546\frac{{\sf
E}X_1^2g(X_1)}{\sigma^2g(\sigma\sqrt{np})}.
$$

\section{The accuracy of the normal approximation to the distributions of Poisson random sums}

In addition to the notation introduced above, let $\lambda>0$ and
$N_{\lambda}$ be the random variable with the Poisson distribution
with parameter $\lambda$:
$$
{\sf P}(N_{\lambda}=k)=e^{-\lambda}\frac{\lambda^k}{k!},\ \ \
k\in\mathbb{N}\cup \{0\}.
$$
Assume that for each $\lambda>0$ the random variables
$N_{\lambda},X_1,X_2,\ldots$ are jointly independent. Consider the
{\it Poisson random sum}
$$
S_{N_{\lambda}}=X_1+\ldots+X_{N_{\lambda}}.
$$
If $N_{\lambda}=0$, then we set $S_{N_{\lambda}}=0$. It is easy to
see that ${\sf E}S_{\lambda}=0$ and ${\sf
D}S_{\lambda}=\lambda\sigma^2$. The accuracy of the normal
approximation to the distributions of Poisson random sum was
considered by many authors, see the historical surveys in
\cite{KorolevShevtsova, ShevtsovaPoisson}. However, the authors are
unaware of any analogs of the Katz--Osipov-type inequalities (1) and
(6) under relaxed moment conditions.

We will obtain a bound for
$$
\Delta_{\lambda}=\sup_x\big|{\sf
P}\big(S_{\lambda}<x\sigma\sqrt{\lambda}\big)-\Phi(x)\big|.
$$
For this purpose fix $\lambda$ and along with $N_{\lambda}$ consider
the random variable $N_{n,p}$ having the binomial distribution with
{\it arbitrary} parameters $n$ and $p\in(0,1]$ such that
$np=\lambda$. As this is so, the reasoning used above implies that
$$
{\sf D}S_{N_{\lambda}}={\sf
D}S_{N_{n,p}}=\sigma^2\lambda=\sigma^2np.
$$
Therefore, by the triangle inequality, in accordance with corollary
6 we have
$$
\Delta_{\lambda}\le\Delta_{n,p}+\sup_x|{\sf
P}(S_{N_{\lambda}}<x)-{\sf P}(S_{N_{n,p}}<x)|\le
$$
$$
\le\frac{1.8546}{\sigma^2}\,{\sf
E}X_1^2\min\bigg\{1,\,\frac{|X_1|}{\sigma\sqrt{np}}\bigg\}+\sup_x\sum_{k=0}^{\infty}{\sf
P}\bigg(\sum_{j=1}^kX_j<x\bigg)\big|{\sf P}(N_{n,p}=k)-{\sf
P}(N_{\lambda}=k)\big|\le
$$
$$
\le\frac{1.8546}{\sigma^2}\,{\sf
E}X_1^2\min\bigg\{1,\,\frac{|X_1|}{\sigma\sqrt{np}}\bigg\}+\sum_{k=0}^{\infty}\big|{\sf
P}(N_{n,p}=k)-{\sf P}(N_{\lambda}=k)\big|.\eqno(24)
$$
Estimate the second term on the right-hand side of (24) by the
Prokhorov inequality \cite{Prokhorov1953} (also see
\cite{Shiryaev1989}, p. 76), according to which
$$
\sum_{k=0}^{\infty}\big|{\sf P}(N_{n,p}=k)-{\sf
P}(N_{\lambda}=k)\big|\le 2p\min\{2,\lambda\},
$$
and obtain that {\it for any} $n$ and $p$ such that $np=\lambda$,
there holds the inequality
$$
\Delta_{\lambda}\le\frac{1.8546}{\sigma^2}\,{\sf
E}X_1^2\min\bigg\{1,\,\frac{|X_1|}{\sigma\sqrt{\lambda}}\bigg\}+2p\min\{2,\lambda\}.\eqno(25)
$$
Now, putting in (25) $p=\lambda/n$ and letting $n\to\infty$, we
obtain the final result:

\smallskip

{\sc Theorem 7}. {\it For any $\lambda>0$}
$$
\Delta_{\lambda}\le\frac{1.8546}{\sigma^2}\,{\sf
E}X_1^2\min\bigg\{1,\,\frac{|X_1|}{\sigma\sqrt{\lambda}}\bigg\}.
$$

\smallskip

Using inequalities (21) -- (23) to estimate $\Delta_{n,p}$ in (24),
we obtain the following result.

\smallskip

{\sc Theorem 8.} {\it Whatever function $g\in\mathcal{G}$ is such
that $\E X_1^2g(X_1)<\infty$, there holds the inequality}
$$
\Delta_{\lambda}\le 1.8546\frac{{\sf
E}X_1^2g(X_1)}{\sigma^2g(\sigma\sqrt{\lambda})}.\eqno(26)
$$

\smallskip

{\sc Remark 2.} The upper bound of the absolute constant used in
theorem 8 is {\it uniform} over the class $\mathcal{G}$. In specific
cases this bound can be considerably sharpened. For example, it is
obvious that $g(x)\equiv|x|\in\mathcal{G}$. For such a function $g$
inequality (26) takes the form of the classical Berry--Esseen
inequality for Poisson random sums, the best current upper bound for
the absolute constant in which is given in \cite{Shevtsova2014}:
$$
\Delta_{\lambda}\le 0.3031\frac{{\sf
E}|X_1|^3}{\sigma^3\sqrt{\lambda}}.\eqno(27)
$$

\section{Convergence rate estimates for mixed Poisson random sums}

\subsection{General results}

In this section we extend the results of the preceding section to
the case where the random number of summands has the mixed Poisson
distribution. For convenience, in this case we introduce an
<<infinitely large>> parameter $n\in\mathbb{N}$ and consider random
variables $N_n^{\star}$ such that for each $n\in\mathbb{N}$
$$
{\sf
P}(N_n^{\star}=k)=\int\limits_{0}^{\infty}e^{-\lambda}\frac{\lambda^k}{k!}d{\sf
P}(\Lambda_n<\lambda),\ \ \ k\in\mathbb{N}\cup\{0\},\eqno(28)
$$
for some positive random variable $\Lambda_n$. For simplicity $n$
may be assumed to be the scale parameter of the distribution of
$\Lambda_n$ so that $\Lambda_n=n\Lambda$ where $\Lambda$ is some
positive <<standard>> random variable in the sense, say, that ${\sf
E}\Lambda=1$ (if the latter exists).

Assume that for each $n\in\mathbb{N}$ the random variable
$N_n^{\star}$ is independent of the sequence $X_1,X_2,\ldots$. As
above, let $S_{N_n^{\star}}=X_1+\ldots+X_{N_n^{\star}}$ and if
$N_n^{\star}=0$, then $S_{N_n^{\star}}=0$.

From (28) it is easily seen that, if ${\sf E}\Lambda_n<\infty$, then
${\sf E}N_n^{\star}={\sf E}\Lambda_n$ so that ${\sf
D}S_n=\sigma^2{\sf E}\Lambda_n$.

Let $N_{\lambda}$ be the random variable with the Poisson
distribution with parameter $\lambda$ independent of
$X_1,X_2,\ldots$ For any $x\in\mathbb{R}$ we have
$$
{\sf P}\big(S_{N_n^{\star}}<x\sigma\sqrt{{\sf
E}\Lambda_n}\big)=\sum_{k=0}^{\infty}{\sf P}(N_n^{\star}=k){\sf
P}\big(S_k<x\sigma\sqrt{{\sf E}\Lambda_n}\big)=
$$
$$
=\sum_{k=0}^{\infty}{\sf P}\big(S_k<x\sigma\sqrt{{\sf
E}\Lambda_n}\big)\int\limits_{0}^{\infty}{\sf P}(N_{\lambda}=k)d{\sf
P}(\Lambda_n<\lambda)=
$$
$$
=\int\limits_{0}^{\infty}{\sf
P}\big(S_{N_{\lambda}}<x\sigma\sqrt{{\sf E}\Lambda_n}\big)d{\sf
P}(\Lambda_n<\lambda)=
\int\limits_{0}^{\infty}{\sf
P}\bigg(\frac{S_{N_{\lambda}}}{\sigma\sqrt{\lambda}}<x\sqrt{\frac{{\sf
E}\Lambda_n}{\lambda}}\bigg)d{\sf P}(\Lambda_n<\lambda)=
$$
$$
=\int\limits_{0}^{\infty}\Phi\bigg(x\sqrt{\frac{{\sf
E}\Lambda_n}{\lambda}}\bigg)d{\sf P}(\Lambda_n<\lambda)+
\int\limits_{0}^{\infty}{\sf
P}\bigg(\frac{S_{N_{\lambda}}}{\sigma\sqrt{\lambda}}<x\sqrt{\frac{{\sf
E}\Lambda_n}{\lambda}}\bigg)d{\sf
P}(\Lambda_n<\lambda)-\int\limits_{0}^{\infty}\Phi\bigg(x\sqrt{\frac{{\sf
E}\Lambda_n}{\lambda}}\bigg)d{\sf P}(\Lambda_n<\lambda).\eqno(29)
$$
From (29) it follows that
$$
\Delta^{\star}_n\equiv\sup_x\bigg|{\sf
P}\big(S_{N_n^{\star}}<x\sigma\sqrt{{\sf
E}\Lambda_n}\big)-\int\limits_{0}^{\infty}
\Phi\Big(\frac{x}{\sqrt{\lambda}}\Big)d{\sf
P}\big(\Lambda_n<\lambda{\sf E}\Lambda_n\big)\bigg|\le
$$
$$
\le\int\limits_{0}^{\infty}\sup_x\bigg|{\sf
P}\bigg(\frac{S_{N_{\lambda}}}{\sigma\sqrt{\lambda}}<x\bigg)-\Phi(x)\bigg|d{\sf
P}(\Lambda_n<\lambda)\le
\int\limits_{0}^{\infty}\Delta_{\lambda}d{\sf
P}(\Lambda_n<\lambda).\eqno(30)
$$

Now, if to estimate the integrand $\Delta_{\lambda}$ in (30) we use
theorem 7 and recall the notation $F(x)={\sf P}(X<x)$, then by the
Fubini theorem we arrive at the representation
$$
\Delta^{\star}_n\le\frac{1.8546}{\sigma^2}\int\limits_{0}^{\infty}\!\!{\sf
E}X_1^2\min\bigg\{1,\frac{|X_1|}{\sigma\sqrt{\lambda}}\bigg\}d{\sf
P}(\Lambda_n<\lambda)=
\frac{1.8546}{\sigma^2}\int\limits_{0}^{\infty}
\bigg[\int\limits_{-\infty}^{\infty}\!\!x^2\min\Big\{1,\frac{|x|}{\sigma\sqrt{\lambda}}\Big\}dF_1(x)\bigg]d{\sf
P}(\Lambda_n<\lambda)=
$$
$$
=\frac{1.8546}{\sigma^2}\int\limits_{-\infty}^{\infty}x^2\bigg[\int\limits_{0}^{\infty}
\min\Big\{1,\,\frac{|x|}{\sigma\sqrt{\lambda}}\Big\}d{\sf
P}(\Lambda_n<\lambda)\bigg]dF_1(x).\eqno(31)
$$
For $x\in\mathbb{R}$ introduce the function
$$
G_n(x)={\sf
E}\min\Big\{1,\,\frac{|x|}{\sigma\sqrt{\Lambda_n}}\Big\}={\sf
P}\Big(\Lambda_n<\frac{x^2}{\sigma^2}\Big)+\frac{|x|}{\sigma}{\sf
E}\frac{1}{\sqrt{\Lambda_n}}\I\Big(\Lambda_n\ge\frac{x^2}{\sigma^2}\Big).\eqno(32)
$$
The expectation in (32) exists since the random variable under the
expectation sign is bounded by 1. Of course, the particular form of
$G_n(x)$ depends on the particular form of the distribution of
$\Lambda_n$. From (30), (31) and (32) we obtain the following
statement.

\smallskip

{\sc Theorem 9.} {\it If ${\sf E}\Lambda_n<\infty$, then}
$$
\Delta^{\star}_n\le\frac{1.8546}{\sigma^2}{\sf
E}X_1^2G_n(X_1)=\frac{1.8546}{\sigma^2}{\sf
E}X_1^2\min\bigg\{1,\,\frac{|X_1|}{\sigma\sqrt{\Lambda_n}}\bigg\}=
$$
$$
=\frac{1.8546}{\sigma^2}\bigg[{\sf
E}X_1^2\I\big(|X_1|\ge\sigma\sqrt{\Lambda_n}\big)+{\sf
E}\frac{|X_1|^3}{\sigma\sqrt{\Lambda_n}}\I\big(|X_1|<\sigma\sqrt{\Lambda_n}\big)\bigg],
$$
{\it where 
the random variables $X_1$ and $\Lambda_n$ are assumed independent.}

\smallskip

In the subsequent sections we will consider special cases where
$\Lambda_n$ has the exponential, gamma and inverse gamma
distributions.

\subsection{Estimates of the rate of convergence of the distributions
of geometric random sums to the Laplace law}

In this section we consider sums of a random number of independent
random variables in which the number of summands $N_n^{\star}$ has
the geometric distribution with parameter $p=\frac{1}{1+n}$,
$n\in\mathbb{N}$:
$$
{\sf P}(N_n^{\star}=k)=\frac{1}{n+1}\Big(\frac{n}{n+1}\Big)^k,\ \ \
k\in\mathbb{N}\cup\{0\}.\eqno(33)
$$
As usual, we assume that for each $n\in\mathbb{N}$ the random
variables $N_n^{\star},X_1,X_2,\ldots$ are independent. We again use
the notation $S_{N_n^{\star}}=X_1+\ldots+X_{N_n^{\star}}$. If $N_n^{\star}=0$, then we set $S_{N_n^{\star}}=0$. It is
easy to see that ${\sf E}N_n^{\star}=n$, ${\sf
D}S_{N_n^{\star}}=n\sigma^2$. Note that for any
$k\in\mathbb{N}\cup\{0\}$
$$
{\sf P}(N_n^{\star}=k)=\frac{1}{n}\int\limits_{0}^{\infty}{\sf
P}(N_{\lambda}=k)\exp\Big\{-\frac{\lambda}{n}\Big\}d\lambda,
$$
where $N_{\lambda}$ is the random variable with the Poisson
distribution with parameter $\lambda$. This means that for
$N_n^{\star}$ representation (28) holds with $\Lambda_n$ being an
exponentially distributed random variable with parameter $\frac1n$.

In what follows we will use traditional notation
$$
\Gamma(\alpha,z)\equiv\int\limits_{z}^{\infty}y^{\alpha-1}e^{-y}dy,\
\ \ \gamma(\alpha,z)\equiv\int\limits_{0}^{z}y^{\alpha-1}e^{-y}dy,\
\ \text{ and } \ \
\Gamma(\alpha)\equiv\Gamma(\alpha,0)=\gamma(\alpha,\infty)
$$
for upper incomplete gamma-function, lower incomplete gamma-function
and gamma-function itself, respectively, where $\alpha>0$, $z>0$.

In the case under consideration
$$
\frac{1}{n}\int\limits_{0}^{\infty}\Phi\bigg(x\sqrt{\frac{n}{\lambda}}\bigg)\exp\Big\{-\frac{\lambda
}{n}\Big\}d\lambda=\int\limits_{0}^{\infty}\Phi\Big(\frac{x}{\sqrt{y}}\Big)e^{-y}dy=\mathcal{L}(x),
$$
where $\mathcal{L}(x)$ is the Laplace distribution function
corresponding to the density
$$
\ell(x)=\frac{1}{\sqrt{2}}e^{-\sqrt{2}|x|},\ \ \ x\in\mathbb{R}
$$
(see, e. g., lemma 12.7.1 in \cite{KorolevBeningShorgin2011}).

At the same time, the function $G_n(x)$ (see (32)) has the form
$$
G_n(x)= 1-\exp\Big\{-\frac{x^2}{n\sigma^2}\Big\}+
\frac{|x|}{n\sigma}\int\limits_{x^2/\sigma^2}^{\infty}\frac{e^{-\lambda/n}}{\sqrt{\lambda}}d\lambda=
\gamma\Big(1,\frac{x^2}{n\sigma^2}\Big)+\frac{|x|}{\sigma\sqrt{n}}\Gamma\Big(\frac12,\,\frac{x^2}{n\sigma^2}\Big).
$$

So, from theorem 9 we obtain the following result.

\smallskip

{\sc Corollary 8.} {\it Let $N_n^{\star}$ have the geometric
distribution $(33)$. Then}
$$
\sup_x\big|{\sf
P}(S_{N_n^{\star}}<x\sigma\sqrt{n})-\mathcal{L}(x)\big|\le
\frac{1.8546}{\sigma^2}\bigg\{{\sf E}\Big[X_1^2
\gamma\Big(1,\frac{X_1^2}{n\sigma^2}\Big)\Big]+\frac{1}{\sigma\sqrt{n}}{\sf
E}\Big[|X_1|^3\Gamma\Big(\frac12,\,\frac{X_1^2}{n\sigma^2}\Big)\Big]\bigg\}.
$$

\subsection{Estimates of the rate of convergence of the distributions
of negative binomial random sums to the variance-gamma law}

The case more general than that considered in the preceding section
is the case of negative binomial random sums.

Let $r>0$ be an arbitrary number. Assume that representation (28)
holds with $\Lambda_n$ being a gamma-distributed random variable
with the density
$$
p(\lambda)=\frac{\lambda^{r-1}e^{-\lambda/n}}{n^r\Gamma(r)}\ \
\lambda>0.
$$
Then the random variable $N_n^{\star}$ has the negative binomial
distribution with parameters $r$ and $\frac{1}{n+1}$:
$$
{\sf
P}(N_n^{\star}=k)=\frac{1}{n^r\Gamma(r)}\int\limits_{0}^{\infty}e^{-\lambda}\frac{\lambda^k}{k!}\lambda^{r-1}e^{-\lambda/n}d\lambda=
\frac{\Gamma(r+k)}{\Gamma(r)\,k!}\Big(\frac{1}{1+n}\Big)^r\Big(\frac{n}{1+n}\Big)^k,\
\ \ \ \ k\in\mathbb{N}\cup\{0\}.\eqno(34)
$$
Let
$$
\mathcal{V}_r(x)\equiv\frac{1}{\Gamma(r)}\int\limits_{0}^{\infty}\!\Phi\Big(\frac{x}{\sqrt{\lambda}}\Big)\lambda^{r-1}e^{-\lambda}d\lambda,\
\ \ x\in\mathbb{R},
$$
be the symmetric variance-gamma distribution with shape parameter
$r$ (see, e. g., \cite{MadanSeneta1990}).

In the case under consideration ${\sf E}N_n^{\star}={\sf
E}\Lambda_n=nr$ so that ${\sf D}S_{N_n^{\star}}=nr\sigma^2$ and for
any $x\in\mathbb{R}$
$$
\int\limits_{0}^{\infty}\!\Phi\Big(x\sqrt{\frac{{\sf
E}\Lambda_n}{\lambda}}\Big)d{\sf
P}(\Lambda_n<\lambda)=\frac{1}{n^r\Gamma(r)}\int\limits_{0}^{\infty}\!\Phi\Big(x\sqrt{\frac{nr}{\lambda}}\Big)
\lambda^{r-1}e^{-\lambda/n}d\lambda=
$$
$$
=\frac{1}{\Gamma(r)}\int\limits_{0}^{\infty}\!\Phi\Big(\frac{x\sqrt{r}}{\sqrt{\lambda}}\Big)\lambda^{r-1}e^{-\lambda}d\lambda\equiv\mathcal{V}_r(x\sqrt{r}).
$$

Here the function $G_n(x)$ (see (32)) has the form
$$
G_n(x)=
=\frac{1}{n^r\Gamma(r)}\int\limits_{0}^{x^2/\sigma^2}\lambda^{r-1}e^{-\lambda/n}d\lambda+
\frac{|x|}{\sigma
n^r\Gamma(r)}\int\limits_{x^2/\sigma^2}^{\infty}\lambda^{r-3/2}e^{-\lambda/n}d\lambda=
$$
$$
=\frac{1}{\Gamma(r)}\Big[\gamma\Big(r,\frac{x^2}{n\sigma^2}\Big)+\frac{|x|}{\sigma\sqrt{n}}\Gamma\Big(r-\frac12,\frac{x^2}{n\sigma^2}\Big)\Big].
$$

So, from theorem 9 we obtain the following result.

\smallskip

{\sc Corollary 9.} {\it Let $N_n^{\star}$ have the negative binomial
distribution $(34)$. Then}
$$
\sup_x\big|{\sf
P}(S_{N_n^{\star}}<x\sigma\sqrt{n})-\mathcal{V}_r(x)\big|\le
\frac{1.8546}{\sigma^2\Gamma(r)}\bigg\{{\sf
E}\Big[X_1^2\gamma\Big(r,\frac{X_1^2}{n\sigma^2}\Big)\Big]+\frac{1}{\sigma\sqrt{n}}{\sf
E}\Big[|X_1|^3\Gamma\Big(r-\frac12,\frac{X_1^2}{n\sigma^2}\Big)\Big]\bigg\}.
$$

\subsection{Estimates of the rate of convergence of the distributions
of Poisson-inverse gamma random sums to the Student distribution}

Let $r>1$ be an arbitrary number. Assume that representation (28)
holds with $\Lambda_n$ being an inverse-gamma-distributed random
variable with parameters $\frac{r}{2}$ and $\frac{n}{2}$ having the
density
$$
p(\lambda)=\frac{n^{r/2}\lambda^{-r/2-1}}{2^{r/2}\Gamma(\frac{r}{2})}\exp\Big\{-\frac{n}{2\lambda}\Big\},\
\  \ \lambda>0.
$$
Then the random variable $N_n^{\star}$ has the so-called
Poisson-inverse gamma distribution:
$$
{\sf
P}(N_n^{\star}=k)=\frac{n^{r/2}}{2^{r/2}\Gamma(\frac{r}{2})}\int\limits_{0}^{\infty}e^{-\lambda}\frac{\lambda^k}{k!}\lambda^{-r/2-1}
\exp\Big\{-\frac{n}{2\lambda}\Big\}d\lambda,\ \ \ \ \
k\in\mathbb{N}\cup\{0\},\eqno(35)
$$
which is a special case of the so-called Sichel distribution see, e.
g., \cite{Sichel1971, Willmot1993}. In this case
$$
{\sf E}\Lambda_n=\frac{n}{r-2}
$$
so that
$$
{\sf D}S_n^{\star}=\frac{n\sigma^2}{r-2}.
$$
Nevertheless, we will normalize random sums not by their mean square
deviations, but by slightly different and asymptotically equivalent
quantities $\sigma\sqrt{n/r}$.

As is known, if $\Lambda_n$ has the inverse gamma distribution with
parameters $\frac{r}{2}$ and $\frac{n}{2}$, then $\Lambda_n^{-1}$
has the gamma distribution with the same parameters. Therefore, we
have
$$
\frac{n^{r/2}}{\Gamma(\frac{r}{2})}\int\limits_{0}^{\infty}\Phi\Big(x\sqrt{\frac{n}{r\lambda}}\Big)\lambda^{-r/2-1}\exp\Big\{-\frac{n}{2\lambda}\Big\}d\lambda=
\frac{n^{r/2}}{\Gamma(\frac{r}{2})}\int\limits_{0}^{\infty}\Phi\Big(x\sqrt{\frac{n\lambda}{r}}\Big)\lambda^{r/2-1}\exp\Big\{-\frac{n\lambda}{2}\Big\}d\lambda=
$$
$$
=\frac{1}{2^{r/2}\Gamma(\frac{r}{2})}\int\limits_{0}^{\infty}\Phi\Big(x\sqrt{\frac{\lambda}{r}}\Big)\lambda^{r/2-1}e^{-\lambda/2}d\lambda=\mathcal{T}_r(x),\
\ \ x\in\mathbb{R},
$$
where $\mathcal{T}_r(x)$ is the Student distribution function with
parameter $r$ ($r$ <<degrees of freedom>>) corresponding to the
density
$$
t_r(x)=\frac{\Gamma(\frac{r+1}{2})}{\sqrt{\pi
r}\Gamma(\frac{r}{2})}\Big(1+\frac{x^2}{r}\Big)^{-(r+1)/2},\ \ \
x\in\mathbb{R},
$$
see, e. g., \cite{BeningKorolev2004}.

In this case the function $G_n(x)$ (see (32)) has the form
$$
G_n(x)={\sf
P}\Big(\Lambda_n^{-1}>\frac{\sigma^2}{x^2}\Big)+\frac{|x|}{\sigma}{\sf
E}\sqrt{\Lambda_n^{-1}}\I\Big(\Lambda_n^{-1}\le\frac{\sigma^2}{x^2}\Big)=
$$
$$
=\frac{n^{r/2}}{2^{r/2}\Gamma(\frac{r}{2})}\int\limits_{\sigma^2/x^2}^{\infty}\lambda^{r/2-1}e^{-n\lambda/2}d\lambda+
 \frac{|x|n^{r/2}}{2^{r/2}\sigma\Gamma(\frac{r}{2})}\int\limits_{0}^{\sigma^2/x^2}\lambda^{(r-1)/2}e^{-n\lambda/2}d\lambda=
$$
$$
=
\frac{1}{\Gamma(\frac{r}{2})}\Big[\Gamma\Big(\frac{r}{2},\frac{n\sigma^2}{2x^2}\Big)+
\frac{|x|}{\sigma}\sqrt{\frac{n}{2}}\,\gamma\Big(\frac{r+1}{2},\frac{n\sigma^2}{2x^2}\Big)\Big],
$$
where $\gamma(\,\cdot\,,\,\cdot\,)$ and
$\Gamma(\,\cdot\,,\,\cdot\,)$ are the lower and upper incomplete
gamma-functions, respectively. So, from theorem 9 we obtain the
following result.

\smallskip

{\sc Corollary 10.} {\it Let $N_n^{\star}$ have the Poisson-inverse
gamma distribution $(35)$. Then}
$$
\Delta_n^{\star}=\sup_x\Big|{\sf
P}\Big(S_{N_n^{\star}}<x\sigma\sqrt{\frac{n}{r}}\Big)-\mathcal{T}_r(x)\Big|\le
\frac{1.8546}{\sigma^2\Gamma(\frac{r}{2})}\bigg\{{\sf
E}\Big[X_1^2\Gamma\Big(\frac{r}{2},\frac{n\sigma^2}{2X_1^2}\Big)\Big]+\frac{1}{\sigma}
\sqrt{\frac{n}{2}}\,{\sf
E}\Big[|X_1|^3\gamma\Big(\frac{r+1}{2},\frac{n\sigma^2}{2X_1^2}\Big)\Big]\bigg\}.
$$

\small

\renewcommand{\refname}{References}


\begin{thebibliography}{99}

\bibitem{BeningKorolev2004} {\it В. Е. Бенинг и В. Ю. Королев.} Об использовании распределения
Стьюдента в задачах теории вероятностей и математической статистики
// Теория вероятностей и ее применения, 2004. Т. 49. Вып.
3. С. 417--435.

\bibitem{Bhat1982} {\it Р.~Н.~Бхаттачария, Р.~Ранга Рао.}
Аппроксимация нормальным распределением. -- М.: Наука, 1982.

\bibitem{Zolotarev1986}
{\it В.\,М. Золотарев.} Современная теория суммирования независимых
случайных величин. -- М.: Наука, 1986.

\bibitem{Kondrik}
{\it А. С. Кондрик, К. В. Михайлов, В. И. Чеботарев.} О равномерной
оценке разности функций распределения / Тезисы докладов XXXI
Дальневосточной школы-семинара им. акад. Е.\,В.\,Золотова,
Владивосток, 2006, С. 16--17.

\bibitem{KorolevBeningShorgin2011} {\it В. Ю. Королев, В. Е. Бенинг, С. Я.
Шоргин.} Математические основы теории риска. 2-е изд., перераб. и
доп. -- М.: ФИЗМАТЛИТ, 2011.

\bibitem{KorolevPopov2011} {\it В. Ю. Королев, С. В. Попов.} Уточнение
оценок скорости сходимости в центральной предельной теореме при
отсутствии моментов порядков, б{\'о}льших второго // Теория
вероятностей и ее применения, 2011. Т. 56. Вып. 4. С. 797--805.

\bibitem{KorolevPopovDAN} {\it В. Ю. Королев, С. В. Попов.} Уточнение оценок
скорости сходимости в центральной предельной теореме при ослабленных
моментных условиях // Доклады Академии наук, 2012. Т. 445. Вып. 3.
C. 265--270.

\bibitem{NefedovaShevtsova2011} {\it Ю. С. Нефедова, И. Г. Шевцова.} О неравномерных
оценках скорости сходимости в центральной предельной теореме // {\it
Теория вероятностей и ее применения}, 2012. Т. 57. Вып. 1. С.
62--97.

\bibitem{Osipov1966} {\it Л. В. Осипов.} Уточнение теоремы Линдеберга //
Теория вероятностей и ее применения, 1966. Т. 11. Вып. 2. С.
339--342.

\bibitem{Petrov1965}
{\it В.~В.~Петров. } Одна оценка отклонения распределения суммы
независимых случайных величин от нормального закона
// Доклады АН СССР, 1965.  Т. 160. Вып. 5. С.
1013--1015.

\bibitem{Petrov1972}
{\it В.~В.~Петров.} Суммы независимых случайных величин. -- М.:
Наука, 1972.

\bibitem{Petrov1987}
{\it В.~В.~Петров.} Предельные теоремы для сумм независимых
случайных величин. -- М.: Наука, 1987.

\bibitem{Prokhorov1953} {\it Ю. В. Прохоров.} Асимптотическое
поведение биномиального распределения // Успехи математических наук,
1953. Т. 8. С. 135--142.

\bibitem{ShevtsovaPoisson} {\it И. Г. Шевцова.} О точности нормальной
аппроксимации для обобщенных пуассоновских распределений // Теория
вероятностей и ее применения, 2013. Т. 58. Вып. 1. С. 152--176.

\bibitem{Shevtsova2013} {\it И. Г. Шевцова.}
Об абсолютных константах в неравенстве Берри--Эссеена и его
структурных и неравномерных уточнениях // Информатика и ее
применения, 2013. Т. 7. Вып. 1. С. 124--125.

\bibitem{Shevtsova2014} {\it И. Г. Шевцова.} Об абсолютных константах в неравенствах типа
Берри-Эссеена // Доклады Академии наук, 2014. Т. 456. Вып. 6. С.
650--654.

\bibitem{Shiryaev1989} {\it А. Н. Ширяев.} Вероятность. -- М.:
Наука, 1989.

\bibitem{BarbourHall1984} {\it A.\,D.\,Barbour, P.\,Hall.} Stein's method
and the Berry--Esseen theorem // Australian Journal of Statistics,
1984. Vol.\,26. P.\,8--15.

\bibitem{ChenShao2001} {\it L. H. Y. Chen, Q. M. Shao.} A non-uniform Berry--Esseen
bound via Stein's method // Probability Theory and Related Fields,
2001. Vol. 120. P. 236--254.

\bibitem{CoxRossRubinstein1979} {\it J. C. Cox, S. A. Ross, M. Rubinstein.} Option pricing:
A simplified approach // Journal of Financial Economics, 1979. Vol.
7. P. 229--263.

\bibitem{Cramer1930} {\it H. Cram{\'e}r.} On the Mathematical Theory of Risk. Skandia
Jubilee Volume, Stockholm Centraltryckeriet, 1930. Reprinted in:
{\it Harald Cramer.} Collective Works. Vol. 1. -- Berlin:
Springer-Verlag, 1994. P. 601--678.

\bibitem{Feller1968}
{\it W. Feller.} On the Berry--Esseen theorem
// Z. Wahrsch. Verw. Geb., 1968.  Bd.\,10. S.\,261--268.

\bibitem{Hoeffding1948}
{\it W. Hoeffding.} The extrema of the expected value of a function
of independent random variables // Ann. Math. Statist., 1948.
Vol.\,19.  P.\,239--325.

\bibitem{Kalashnikov1997} {\it V. V. Kalashnikov.} Geometric Sums:
Bounds for Rare Events with Applications. -- Dordrecht: Kluwer
Academic Publishers, 1997.

\bibitem{Katz1963}
{\it M. Katz.} Note on the Berry--Esseen theorem
// Annals of Math. Statist., 1963. Vol. 39. No. 4.  P. 1348--1349.

\bibitem{KP2011_3}
{\it V. Korolev, S. Popov.} On the universal constant in the
Katz--Petrov and Osipov inequalities // Discussiones Mathematicae.
Probability and Statistics, 2011. Vol.~31. P.~29--39.

\bibitem{KorolevShevtsova} {\it V. Yu. Korolev, I. G. Shevtsova.} An improvement
of the Berry--Esseen inequality with applications to Poisson and
mixed Poisson random sums // Scandinavian Actuarial Journal, 2012.
Vol. 2012. Issue 2. P. 81--105.

\bibitem{MadanSeneta1990} {\it D. B. Madan, E. Seneta}. The variance gamma $($V.G.$)$
model for share market return // Journal of Business, 1990. Vol. 63.
P. 511--524.

\bibitem{Paditz1980} {\it L. Paditz.} Bemerkungen zu einer
Fehlerabsch{\"a}tzung im zentralen Grenzwertsatz // Wiss. Z.
Hochschule f{\"u}r Verkehrswesen <<Friedrich List>>, 1980. Bd. 27.
\textnumero\, 4. S. 829--837.

\bibitem{Paditz1984} {\it L. Paditz.} On error-estimates in the central limit
theorem for generalized linear discounting // Math.
Operationsforsch. u. Statist., Ser. Statistics, 1984. Bd. 15. \textnumero\,
4. S. 601--610.

\bibitem{Paditz1986} {\it L. Paditz.} {\"U}ber eine Fehlerabsch{\"a}tzung im
zentralen Grenzwertsatz // Wiss. Z. Hochschule f{\"u}r Verkehrswesen
<<Friedrich List>>, 1986. Vol. 33. No. 2. P. 399--404.

\bibitem{Sichel1971} {\it H. S. Sichel.} On a family of discrete distributions
particular suited to represent long tailed frequency data /
N.~F.~Laubscher (Ed.). Proceedings of the 3rd Symposium on
Mathematical Statistics. -- Pretoria: CSIR, 1971. P. 51--97.

\bibitem{Sunklodas} {\it J. K. Sunklodas.} On the normal approximation of a
binomial random sum // Lithuanian Mathematical Journal, 2014. Vol.
54. No. 3. DOI:10.1007/s10986-014-9248-6.

\bibitem{Willmot1993} {\it G. E. Willmot.} On recursive evaluation of mixed
Poisson probabilities and related quantities // Scandinavian
Actuarial Journal, 1993, No. 2. P. 114--133.

\bibitem{Zolotarev1997}
{\it V. M. Zolotarev.} Modern Theory of Summation of Random
Variables. -- Utrecht: VSP, 1997.

\end{thebibliography}
\end{document}